\documentclass[12pt]{article}
\usepackage[utf8]{inputenc}
\usepackage{cite}
\usepackage{enumerate}
\usepackage{graphicx}
\usepackage{amsmath,amssymb}
\usepackage{color}
\usepackage{hyperref}
\usepackage{authblk}

\title{On the 4-color theorem for signed graphs}
\author[1]{Franti\v sek Kardo\v s \footnote{This research was supported by the French National Research Agency (ANR) project HOSIGRA (HOmomorphisms of SIgned GRAphs) no. ANR-17-CE40-0022.}}
\author[1]{Jonathan Narboni}
\affil[1]{Univ. Bordeaux, CNRS,  Bordeaux INP, LaBRI, UMR 5800, F-33400, Talence, France}
\date{\today}
\newtheorem{theorem}{Theorem}
\newtheorem{definition}{Definition}
\newtheorem{conjecture}{Conjecture}
\newtheorem{corollary}{Corollary}
\newtheorem{lemma}{Lemma}

\newcommand{\pf}{{{\it Proof.}\ }}
\newcommand{\ep}{\hfill$\square$}
\begin{document}

\maketitle

\begin{abstract}
There are several ways to generalize graph coloring to signed graphs. M\'a\v cajov\'a, Raspaud and \v Skoviera introduced one of them and conjectured that in this setting, for signed planar graphs four colours are always enough, generalising thereby The Four Color Theorem. We disprove the conjecture.
\end{abstract}
\section{Introduction}

All graphs in this paper are finite, simple undirected graphs. Let $G$ be a graph and $\sigma : E(G)\rightarrow \{-1,+1\}$ be a mapping. The pair $(G,\sigma)$ is called a \emph{signed graph}, $\sigma$ is called the \emph{signature} of the graph, and $G$ the \emph{underlying graph}. For convenience, when there is no ambiguities, we only write $G$ to denote the signed graph.

Signed graph is a notion introduced by Zaslavsky \cite{zaslavsky1982signed}. He also defined a $k$-signed coloring of a signed graph $G$ (see \cite{Z82}) as a mapping $$c:V(G)\rightarrow\{-k,-(k-1)...-1,0,1...,k-1,k\}$$ such that for every edge $e=uv\in E(G)$, $c(u)\neq\sigma(e)\cdot c(v)$. 

This definition is relevant and natural because it is compatible with the \emph{switch} operation on a signed graph. Switching a vertex $v$ in a signed graph consists in switching the sign of every edge incident to $v$. This operation induces equivalence classes on signed graphs that have the same underlying graph.

Because switching a vertex preserves the sign of every cycle in the graph (i.e., the product of the signs of the edges composing the cycle), the equivalence switching classes can be characterized by the sign of the cycles in the graphs: Two signed graphs $G_1$ and $G_2$ are equivalent if and only if every cycle of $G_1$ have the same sign as the corresponding cycle of $G_2$. 

A signed coloring corresponding to the Zaslavsky's definition is also preserved by the switch operation. When switching a vertex $v$ in a colored signed graph $G$, switching the sign of the color of $v$ preserves the coloring of $G$, but with this definition, a $k$-signed coloring of $G$ uses $2k+1$ colors.

M\'a\v cajov\'a, Raspaud and \v Skoviera \cite{MRS} introduced a new definition of $k$-signed coloring in which only $k$ colors are used. According to this definition, a $k$-signed coloring of a signed graph $G$ is, if $k\equiv 0[2]$ (resp. $k\equiv 1[2]$) a mapping $c:V(G)\rightarrow\{-k/2,-(k/2-1),...,-1,1,...,k/2-1,k/2\}$ (resp. $c:V(G)\rightarrow\{-(k-1)/2,-((k-1)/2-1),...,-1,0,1,...,(k-1)/2-1,(k-1)/2\}$), such that for every edge $e=uv\in E(G)$, $c(u)\neq\sigma(e)\cdot c(v)$. 

Given a signed graph $(G,\sigma)$, we denote by $\chi(G,\sigma)$ (or $\chi(G)$) the chromatic number of $G$: $\chi(G)=min\{k:\ $$G$ has a $k$-signed coloring$\}$.

M\'a\v cajov\'a, Raspaud and \v Skoviera also conjectured that the Four Color Theorem holds for the signed planar graphs as well :

\begin{conjecture} \cite{MRS}
Let $G$ be a simple signed planar graph. Then $\chi(G)\le 4.$
\label{conj:raspaud}
\end{conjecture}

Signed coloring is also closely related to list-coloring, and Conjecture \ref{conj:raspaud} would in fact imply another conjecture 
about a special type of list-coloring of (non-signed) graphs called \textit{weak list-coloring}.

A list assignement $L$ of a graph $G$ is \textit{symmetric} if,  for every vertex $v$ of $G$, $L(v)$ is such that, if an integer $i$ is in $L(v)$, then $-i$ is also in $L(v)$.
A graph $G$ is \textit{$k$-weakly choosable} if, for every symmetric list-assignement $L$, such that for every vertex $v$, $|L(v)|\leq k$, then $G$ is $L$-colorable. The \textit{weak choice number} of a graph $G$ is the minimum $k$ such that $G$ is $k$-weakly choosable, this number is denoted $ch^w(G)$.

 Zhu
 \cite{zhu2018refinement} proposed the following generalization of the Four Color Theorem:

\begin{conjecture} \cite{zhu2018refinement}
Let $G$ be a planar graph. Then $ch^w(G)\leq 4$.

\label{conj:zhuxuding}
\end{conjecture}

Zhu \cite{zhu2018refinement} 
showed that Conjecture \ref{conj:raspaud} implies Conjecture \ref{conj:zhuxuding}; in this paper, we prove that Conjecture \ref{conj:raspaud} does not hold.

\section{Results}

Before introducing a counterexample to Conjecture \ref{conj:raspaud}, we translate the question of vertex-coloring of a signed planar graph to a question of edge-labeling of its dual. We extend this way the reduction from $4$-coloring of a triangulation to $3$-edge-coloring of the dual used in the proof of the four-color theorem (see \cite{robertson1997four}).
\bigskip

Let $G$ be a 3-connected signed planar graph, and let $f$ be a face of $G$. We call the face $f$ positive (negative) if the facial cycle of $f$ contains an even (odd, respectively) number of negative edges.

Let $G^*$ be the dual graph of $G$. We call the vertices of $G^*$ corresponding to the positive faces of $G$ \emph{positive}, and the vertices corresponding to the negative face of $G$ \emph{negative}.

Observe that there is always an even number of negative faces, and that the set of negative faces is invariant with respect to switching, so signature of the vertices of the dual graph $G^*$ is the same for every graph $G$ belonging to a same switch equivalence class.

In the figures, the positive vertices will be represented by simple dots, and the negative ones will be represented by circles with a minus sign inside.

Let $H$ be a 3-connected planar graph and let $c$ be a $\{0,a,b\}$-edge-labeling of $H$. We denote $d_x(v)$ the number of edges incident to $v$ labelled $x$ for $v\in V(H)$ and $x\in\{0,a,b\}$.

\begin{definition}
Let $H$ be a 3-connected planar graph with an even number of negative vertices and let $c$ be a $\{0,a,b\}$-edge-labeling of $H$. The labeling $c$ is a \emph{weak signed edge-labeling} of $H$ if
\begin{enumerate}
    \item[(i)] $d_0(v) \equiv d_H(v) \pmod{2}$, and
    \item[(ii.a)] $d_a(v) \equiv d_b(v) \equiv d_H(v) \pmod{2}$ if $v$ is a positive vertex, or
    \item[(ii.b)] $d_a(v) \equiv d_b(v) \equiv  d_H(v)+1 \pmod{2}$ if $v$ is a negative vertex.
\end{enumerate}
\end{definition}

In particular, if $H$ is a cubic graph with an even number of negative vertices, then a weak signed edge-labeling of $H$ is a $\{0,a,b\}$-edge-labeling of $H$ such that 
\begin{itemize}
    \item if $v$ is positive, then it is incident to one edge of each label from $\{0,a,b\}$,
    \item if $v$ is negative, then it is incident to one edge labelled $0$, and the two other edges have the same label $c \in \{0,a,b\}$.
\end{itemize}

\begin{conjecture}
Let $H$ be a 3-connected planar graph with an even number of negative vertices. Then there exists a weak signed edge-labeling of $H$.
\label{conj:weak}
\end{conjecture}

\begin{theorem}
Conjectures \ref{conj:raspaud} and \ref{conj:weak} are equivalent.
\label{th:conj_equiv}
\end{theorem}
\pf We prove first that Conjecture \ref{conj:raspaud} implies Conjecture \ref{conj:weak}. Suppose the former true.

Let $H$ be a 3-connected planar graph with an even number of negative vertices. Let $T$ be the set of negative vertices of $H$. Then there exists a $T$-join $F$ (a subgraph of $H$ such that $d_F(v)$ is odd if and only if $v\in T$) in $H$ -- it suffices to consider a binary sum of (not necessarily disjoint) paths joining disjoint pairs of vertices in $T$. 

We will call an edge $e$ of $H$ \emph{negative} if $e\in F$, otherwise it is \emph{positive}. We denote $d_x^+(v)$ ($d_x^-(v)$) the number of positive (negative) edges labeled $x$ incident to $f$, respectively. 

Let $\sigma$ be a signature of $G = H ^*$ defined by
$$
\sigma(e)=\begin{cases} 1 & \textrm{if $e^*\notin F$,} \\ -1 & \textrm{if $e^* \in F$}.
\end{cases}
$$
By definition, a face $f$ of $G$ is negative if and only if the corresponding vertex $f^*$ of $H$ is negative.

Let $\varphi$ be a signed 4-coloring of $(G,\sigma)$ with colors from $\{-2,-1,1,2\}$. 
Let $e=uv$ be an edge of $G$ and $e^*$ be the edge of $G^*$ corresponding to $e$. The label $\varphi^*(e^*)$ of $e^*$ is defined depending on the sign of $e$ and the colors of $u$ and $v$ in the following way:
	\begin{equation}
	\varphi^*(e^*)=
	\begin{cases}
	    0 & \textrm{if $\varphi(u)=-\sigma(uv)\cdot \varphi(v)$,} \\
	    a & \textrm{if $\sigma(uv)\cdot\varphi(u)\cdot\varphi(v)=2$,} \\
	    b & \textrm{if $\sigma(uv)\cdot\varphi(u)\cdot\varphi(v)=-2$.} 
	\end{cases}
	\label{eq:def}
	\end{equation}

Observe that if $\varphi(u) \ne \pm \sigma(uv)\cdot \varphi(v)$, then $\{|\varphi(u)|,|\varphi(v)|\}=\{1,2\},$ and so  $\varphi^*$ is well-defined. It suffices to prove that $\varphi^*$ is a weak signed edge-labeling.
	

Let $e=uv$ be an edge of $G$. When passing from $\varphi(u)$ to $\varphi(v)$, the color may (or may not) change the sign and/or change the absolute value. 

Let $f$  be a face of $G$. Consider the edges of the cycle defining the boundary of $f$.
The following observations are direct consequences of the definition of $\varphi^*$:

\begin{enumerate}
\item Each change of the absolute value of vertex color corresponds to an $a$- or $b$-edge; the number of such changes around $f$ is even. Therefore, there is an even number of $a$- or $b$-labeled edges incident to $f^*$, which is equivalent to $$d_0(f^*)\equiv d_{H}(f^*) \pmod{2},$$ so $\varphi^*$ satisfies the condition $(i)$.
In other words, $$d_a(f^*) \equiv d_b(f^*) \pmod{2}.$$ 
    In particular, if $f$ is a triangle then there is an odd number of $0$-labeled edges incident to $f^*$.
\item Each change of the sign of vertex color corresponds either to a positive $0$- or $b$-labeled edge, or to a negative $a$-labeled edge; the number of such changes around $f$ is even again. Therefore,
    $$ d_0^+(f^*)+d_b^+(f^*)+d_a^-(f^*) \equiv 0 \pmod{2},$$
    which is equivalent to
    $$ d_0^+(f^*)+d_b^+(f^*) +d_a^+(f^*) \equiv d_a^-(f^*) + d_a^+(f^*)  \pmod{2},
    $$
    meaning
    $$
    d^+(f^*)\equiv d_a(f^*) \pmod{2},
    $$
        and so the total number of positive edges incident to $f$ in $G$ has the same parity as the total number of $a$-labeled edges incident to $f^*$ in $G^*$.
        Hence, the labeling $\varphi^*$ satisfies also the condition $(ii)$.
        
        In particular, if $f$ is a positive triangle, then there is an odd number of $a$-labeled edges (and of $b$-labeled edges) incident to $f^*$; if $f$ is a negative triangle, then there is an even number of $a$-labeled edges (and of $b$-labeled edges) incident to $f^*$.
\end{enumerate}


Conversely, let $G,\sigma$ be a 3-connected signed planar graph.
Let $\varphi^*$ be a weak signed $\{0,a,b\}$-edge-labeling of $G^*$. We define a coloring of $G$ in the following way: Let $T$ be a spanning tree of $G$ rooted at a vertex $r$. We set $\varphi(r)=1$ and for every other vertex $u$ of $G$, the color of $u$ is defined depending on the color of its father $v$ in the spanning tree, the label of the edge $e^*$ corresponding to to the edge $e=uv$ and the sign of $e$ in the following way:

\begin{equation}
	\varphi(u)=
	\begin{cases}
	    \sigma(e)\cdot(\sigma(\varphi(v))\cdot3-\varphi(v)) & \textrm{if $\varphi^*(e^*) = a$,} \\
	    -\sigma(e)\cdot(\sigma(\varphi(v))\cdot3-\varphi(v)) & \textrm{if $\varphi^*(e^*) = b$,} \\
	    -\sigma(e)\cdot\varphi(v) & \textrm{if $\varphi^*(e^*) = 0$.} \\
	\end{cases}
\label{eq:dual}
\end{equation}
    
We claim that the coloring $\varphi$ of $G$ defined this way and the coloring $\varphi^*$ of $G^*$ satisfy the equation (\ref{eq:def}) for every edge $e$ of $G$, and so $\varphi$ is a proper signed 4-coloring of $G$.

It is straightforward to verify that the formulae (\ref{eq:def}) and (\ref{eq:dual}) are equivalent for every edge $e\in T$. 

The edges of $G^*$ corresponding to the edges from $E(G)\setminus E(T)$ form a spanning tree of $G^*$. We will prove that (\ref{eq:def}) is verified for these edges by induction.

Let $f^*$ be a vertex of $G^*$ corresponding to a face $f$ of $G$ such that (\ref{eq:def}) has already been verified for all the incident edges but one; let that edge be $e=uv$. 

Since, for each $c\in\{0,a,b\}$, the parity of the $c$-labeled edges incident to $f^*$ is determined by the degree and the sign of $f^*$, the label of $e^*$ is uniquely determined by the labels of the other edges incident to $f^*$ and the sign of $e$.

Observe that $\varphi^*(e^*)\in\{a,b\}$ if and only if there is a change of absolute value of colors of vertices ; at $f^*$ there is always an even number of such edges by $(i)$.

Similarly, an edge represents a change of the sign of colors of vertices if and only if it is a positive $0$- or $b$-labeled edge or a negative $a$-labeled edge; it follows from the definition of a weak signed edge-labeling that there is always an even number of such edges at $f^*$.

Let $P^*$ be the set of edges incident to $f^*$ distinct from $e^*$, let $P$ be the corresponding path from $u$ to $v$ along $f$ in $G$.

By definition of a weak signed edge-labeling, $\varphi^*(e^*)=0$ if and only if among the edges on $P^*$ there is an even number of edges labeled $a$ or $b$, which means that along $P$, there is an even number of changes of absolute value, and so $\varphi(u)$ and $\varphi(v)$ have the same absolute value.

Moreover, $\varphi^*(e^*)$ has a label that corresponds to a change of sign (a positive $0$- or $a$-labeled edge or a negative $b$-labeled edge) if and only if among the edges of $P^*$ there is an odd number of such edges, which means that along $P$, there is an odd number of changes of the sign, and so $\varphi(u)$ and $\varphi(v)$ have different signs.

The last two paragraphs combined together imply that (\ref{eq:def}) is true for the edge $e=uv$.
\ep

\begin{definition}
Let $H$ be a 3-connected planar graph with an even number of negative vertices. A $\{0,a,b\}$-edge-labeling $c$ of $H$ is a \emph{strong signed edge-labeling} if 
\begin{enumerate}
    \item[(i)] $c$ is a weak signed edge-labeling of $H$, and
    \item[(ii)] $d_0(v) < d_H(v)$ for every odd-degree vertex $v$ of $H$.
\end{enumerate}
\end{definition}

Observe that $d_0(v) = d_H(v)$ is possible only if $v$ is a negative vertex of odd degree.

\begin{conjecture}
Let $H$ be a 3-connected planar graph with an even number of negative vertices. Then there exists a strong signed edge-labeling of $H$.
\label{conj:strong_labeling}
\end{conjecture}

\begin{theorem}
    Conjectures \ref{conj:weak} and \ref{conj:strong_labeling} are equivalent.
    \label{th:equiv_weak_strong}
\end{theorem}

\pf
Trivially, Conjecture \ref{conj:strong_labeling} implies Conjecture \ref{conj:weak}.

\begin{figure}
    \centering
    \includegraphics[scale=1]{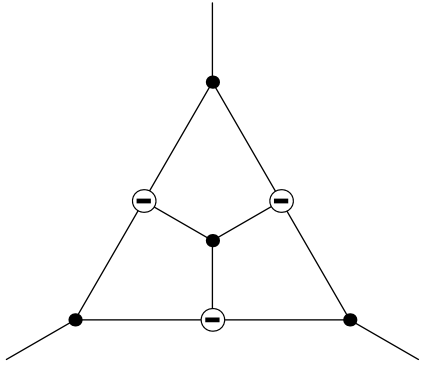}
    \caption{An example of a gadget used to replace an odd negative vertex (here $k=3$).}
    \label{fig:gadget}
\end{figure}

Let $k\ge 3$ be an odd integer. We define $W_k$ as the graph obtained from a cycle of length $2k$ with every other vertex negative by adding a pending edge to each positive vertex of the cycle, and by adding a positive vertex adjacent to each negative vertex of the cycle. See Figure \ref{fig:gadget} for illustration.

We will call the edges of the $2k$-cycle in $W_k$ \emph{ring} edges, and the edges joining negative vertices of the cycle to the central positive vertex \emph{spikes}.

Suppose Conjecture \ref{conj:weak} is true.
Let $H$ be a 3-connected planar graph with an even number of negative vertices. Let $H'$ be the graph obtained from $H$ by replacing every odd negative vertex $v$ by a copy of $W_{d(v)}$. 
Since the graph $H'$ is planar, 3-connected and has an even number of negative vertices, it has a weak signed edge-labeling.

We claim that we can reduce this weak signed edge-labeling of $H'$ to a strong signed edge-labeling of $H$ simply by a contraction of each gadget to a single (negative) vertex, keeping the labels of the edges. 

Let $W$ be a copy of $W_k$ in $H'$ for some odd $k\ge 3$. 

Observe first that not all the edges leaving $W$ are labeled $0$. If this was the case, all the ring edges would be labeled $a$ or $b$. However, along the ring, a positive 3-vertex corresponds to a change from $a$ to $b$ (or vice versa), whereas a negative 3-vertex cannot be incident to an $a$- and $b$-labeled edge at the same time. There has to be an even number of changes, but there is an odd number of positive ring vertices, a contradiction.

We need yet to prove that the number of edges leaving $W$ 
labeled $a$ (resp. $b$) is even; the fact that the number of edges leaving $W$ labeled $0$ is odd will follow automatically.

Let us count the number of incidences with $a$-labeled edges.
Since the central vertex is a positive odd vertex, it is incident to an odd number of $a$-labeled edges. 
There is an odd number of negative vertices on the ring, each of them is incident to an even number of $a$-labeled edges. 
There is an odd number of positive vertices on the ring, each of them is incident to an odd number of $a$-labeled edges. 

Therefore, in total, there is an even number of $a$-labeled edge incidences in $W$, and so there is an even number of edges leaving $W$ labeled $a$.
\ep

As a direct consequence of the definition of the strong signed edge-labeling we get the following: 

\begin{corollary}
Let $H$ be a cubic 3-connected planar graph with an even number of negative vertices. $H$ has a strong signed edge-labeling if and only if $H$ has a 2-factor $F$ such that every cycle of $F$ has an even number of positive vertices. 
\end{corollary}

If $H$ admits such a 2-factor $F$, we call $F$ a \emph{consistent} 2-factor. 

\pf If $H$ has a strong signed edge-labeling, then for each $v\in V(H)$, $d_a(v)+d_b(v) = 2$. Hence, the edges labeled $a$ or $b$ form a 2-factor $F$ of $H$. Moreover, as for $v\in V(H)$, $v$ is a positive vertex if and only if $d_a(v)=d_b(v)=1$, each cycle of $F$ must have an even number of positive vertex.

Conversely, assume that $H$ has a consistent 2-factor $F$. For each cycle $C\in F$, choose an edge $e \in C$, and label $e$ with $a$. Then, label the other edges of $C$ with $a$ and $b$, in such a way that the labels change only at positive vertices. The edges that are not part of $F$ are labeled $0$. It is easy to see that such a labeling is a strong signed edge-labeling of $H$.
\ep

Observe that if $H$ is a hamiltonian cubic planar graph with an even number of negative vertices, then any Hamilton cycle of $H$ is a consistent 2-factor.

\begin{figure}[ht]
    \centering
    \includegraphics[scale=1]{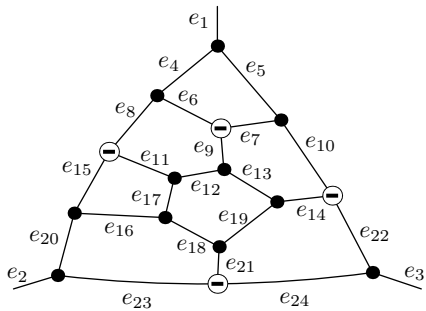}
    \caption{Tutte's fragment}
    \label{fig:fragment}
\end{figure}

\begin{lemma}
    Let $H$ be a cubic 3-connected planar graph with an even number of negative vertices, containing a Tutte's fragment $T_0$ attached by the edges $e_1$, $e_2$, $e_3$, as depicted in Figure \ref{fig:fragment}. Then every consistent 2-factor of $H$ contains the edge $e_1$.
\label{lemma:fragment}
\end{lemma}

\pf
Assume that $e_1$ is not in of $F$. Then $e_4$, $e_5$ are in $F$.
Moreover, as there is an odd number of positive vertices in the fragment, the edges $e_{2}$ and $e_{3}$ have to be in $F$. 

We introduce a sequence of claims, each one being easy to check.
\begin{itemize}
    \item  $e_{21} \in F$. If not, then $e_{23},e_{24}\in F$, and so there would be an odd number of positive vertices left in the fragment to be covered by $F$.
    \item $e_9\in F$. If not, then there would be a 4-cycle in $F$ with three positive vertices.
    \item $e_{17}\in F$. If not, then $e_{11}$, $e_{12}$, $e_{16}$, $e_{18} \in F$, and so $e_{13},e_{19} \notin F$, so $F$ is not a 2-factor.
    \item $e_{15}\in F$. If not, then $e_{8}$, $e_{11}$, $e_{16}$, $e_{20}\in F$, and so $e_6, e_{12}, e_{18}, e_{23} \notin F$, meaning $e_7, e_{13}, e_{19}, e_{24} \in F$, so $F$ does not cover all the vertices of the fragment.
    \item $e_{14}\in F$. If not, then $e_{10}$, $e_{13}$, $e_{19}$, $e_{22}\in F$, and so $e_7, e_{12}, e_{18}, e_{25} \notin F$, meaning $e_{11}, e_{16} \in F$, so $F$ contains a 4-cycle with three positive vertices.
\end{itemize}   

\begin{figure}
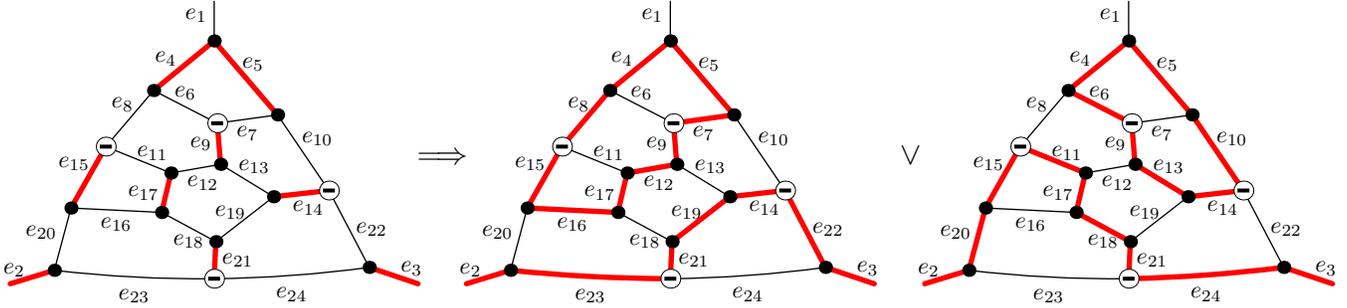

\centerline{
\begin{tabular}{ccccc}
\begin{tabular}{c}
    \includegraphics{tutte_frag/tutte_frag-2.mps}
    \hspace*{-25pt}
\end{tabular} & $\Longrightarrow$ &
\begin{tabular}{c}
    \hspace*{-25pt}
    \includegraphics{tutte_frag/tutte_frag-3.mps}
    \hspace*{-15pt}
\end{tabular} & $\vee$
\begin{tabular}{c}
    \hspace*{-15pt}
    \includegraphics{tutte_frag/tutte_frag-4.mps}
\end{tabular}
\end{tabular}
}
    \caption{In a Tutte's fragment with the given position of four negative vertices, for any 2-factor avoiding the edge $e_1$ there is always at least one cycle containing an odd number of positive vertices.}
    \label{fig:tenth}
\end{figure}

From the previous claims, we have that if $H$ has a consistent 2-factor $F$ s.t.~$e_1\notin F$, then $e_2$, $e_3$, $e_4$, $e_5$, $e_9$, $e_{14}$, $e_{15}$, $e_{17}$ and $e_{21}$ are in $F$, as depicted in Figure \ref{fig:tenth}, left. The remaining edges form a cycle, so we only have two choices to complete $F$. Each of these leads to a cycle with an odd number of positive vertices (see Figure \ref{fig:tenth}).  
\ep

\begin{figure}
    \centering
    \includegraphics[scale=1]{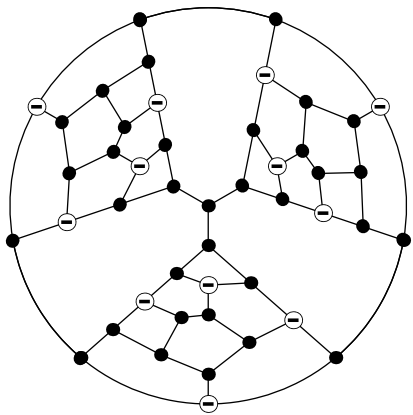}
    \caption{The Tutte's graph composed of three Tutte's fragments.}
    \label{fig:first}
\end{figure}

\begin{theorem}
There exists a set of twelve vertices of the Tutte's graph $T$ to be chosen negative such that $T$ does not have a consistent 2-factor.
    \label{th:counterexample}
\end{theorem}
\pf 
Let the negative vertices of $T$ be chosen as in Figure \ref{fig:first}. 
Assume that $H$ has a consistent 2-factor $F$. 

The graph $T$ can be viewed as a $K_4$ where three of the four vertices were replaced by the Tutte's fragment. By Lemma \ref{lemma:fragment}, all the three edges incident to the central vertex belong to $F$, a contradiction. \ep

\begin{corollary}
Conjecture \ref{conj:raspaud} is false.
\end{corollary}

To find a counterexample, it suffices to consider the Tutte's graph $T$ with a choice of negative vertices as depicted in Figure \ref{fig:first} and replace every negative vertex by the graph $W_3$ depicted in Figure \ref{fig:gadget}, and then take the dual. This gives a graph on 61 vertices.

\section{Concluding remarks}
\hspace*{\parindent}The question that naturally arises is the size of a minimum non-$4$-colorable signed planar graph, and this question remains open. Clearly, it suffices to search for a triangulation whose dual is a non-hamiltonian 3-connected cubic planar graph, and then search for a position of an even number of negative vertices such that there is no weak edge-labeling of the dual graph.

It is known \cite{HM88} that 3-connected cubic planar graphs on at most 36 vertices are all hamiltonian. There are six smallest non-hamiltonian 3-connected cubic planar graphs on 38 vertices, and for each of them it is possible to choose a position of eight negative vertices such that the graph does not admit a consistent 2-factor, and therefore it that does not admit a strong edge-labeling. (We omit the details). 

\begin{figure}
    \centering
    \includegraphics[scale=1]{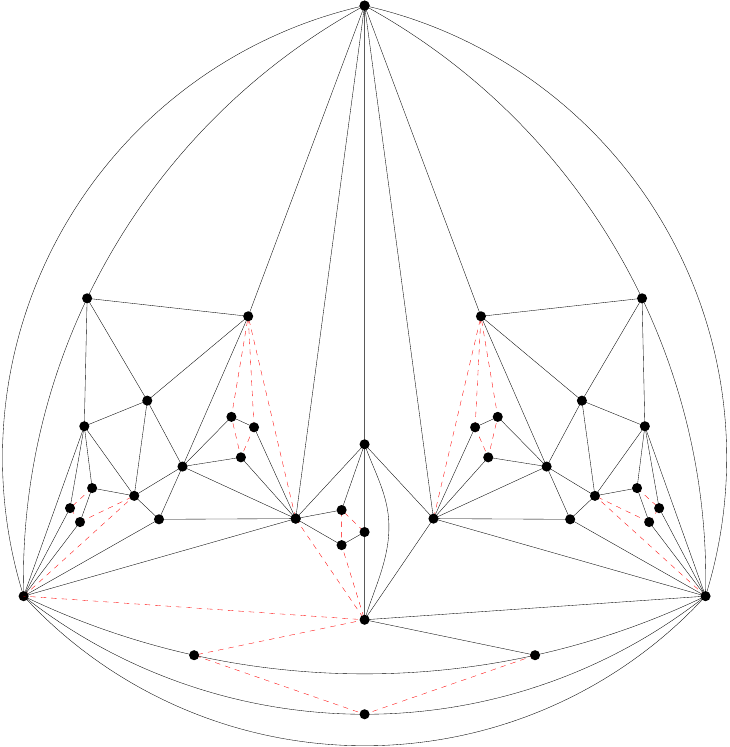}
    \caption{The smallest non-4-colorable signed planar graph we have found so far (dashed red lines stand for negative edges).}
    \label{fig:smallest}
\end{figure}

To guarantee the non-existence of a weak edge-labeling, it suffices to replace four negative vertices by a gadget that has 7 vertices, the corresponding graph that does not admit a weak edge-labeling then has 74 vertices, which corresponds to a non-4-colorable signed triangulation on 39 vertices, see Figure \ref{fig:smallest}. (Again, we omit the details).

Another interesting question is the complexity of deciding whether or not a planar signed graph is $4$-signed colorable. By the way, Conjecture \ref{conj:zhuxuding} remains a challenging and interesting open question.

\bibliographystyle{alpha}
\bibliography{biblio}

\begin{thebibliography}{RSST97}

\bibitem[HM88]{HM88}
Derek~A Holton and Brendan~D McKay.
\newblock The smallest non-hamiltonian 3-connected cubic planar graphs have 38
  vertices.
\newblock {\em Journal of Combinatorial Theory, Series B}, 45(3):305--319,
  1988.

\bibitem[MR{\v{S}}16]{MRS}
Edita M{\'a}{\v{c}}ajov{\'a}, Andr{\'e} Raspaud, and Martin {\v{S}}koviera.
\newblock The chromatic number of a signed graph.
\newblock {\em The Electronic Journal of Combinatorics}, 23(1):P1--14, 2016.

\bibitem[RSST97]{robertson1997four}
Neil Robertson, Daniel Sanders, Paul Seymour, and Robin Thomas.
\newblock The four-colour theorem.
\newblock {\em Journal of combinatorial theory, Series B}, 70(1):2--44, 1997.

\bibitem[Zas82a]{Z82}
Thomas Zaslavsky.
\newblock Signed graph coloring.
\newblock {\em Discrete Mathematics}, 39(2):215--228, 1982.

\bibitem[Zas82b]{zaslavsky1982signed}
Thomas Zaslavsky.
\newblock Signed graphs.
\newblock {\em Discrete Applied Mathematics}, 4(1):47--74, 1982.

\bibitem[Zhu18]{zhu2018refinement}
Xuding Zhu.
\newblock A refinement of choosability of graphs.
\newblock {\em arXiv preprint arXiv:1811.08587}, 2018.

\end{thebibliography}

\end{document}